\documentclass[12pt]{article}
\usepackage{amsmath,amsthm,amsfonts,amssymb,amscd,abstract,titling,enumitem,latexsym,mathrsfs,parskip,geometry,authblk,url}
\newtheorem{lemma}{Lemma}
\newtheorem{remark}{Remark}

\newtheorem{fact}{Fact}
\newtheorem{theorem}{Theorem}
\newtheorem{corollary}{Corollary}

\def\endproof{$\hfill \blacksquare$}
	
\title{Locally nilpotent polynomials over $\mathbb{Z}$}
\author{Sayak Sengupta\\
        Department of Mathematics and Statistics,\\
		Binghamton University - SUNY,\\
		Binghamton, New York 13902-6000,\\
		USA,\\
		\texttt{sengupta@math.binghamton.edu}}
\date{}

\begin{document}
\maketitle

\begin{abstract}
   For a polynomial $u(x)$ in $\mathbb{Z}[x]$ and $r\in\mathbb{Z}$, we consider the orbit of $u(x)$ at $r$, $\mathcal{O}_u(r):=\{u(r),u(u(r)),\ldots\}$. We ask two questions here: (i) what are the polynomials $u$ for which $0\in \mathcal{O}_u(r)$ and (ii) what are the polynomials for which $0\not\in \mathcal{O}_u(r)$ but, modulo every prime $p$, $0\in \mathcal{O}_u(r)$? In this paper we classify the polynomials for which (ii) holds.
   We also present some results for some special $r'$s for which (i) can be answered.
\end{abstract}



\section{Introduction}
A. Borisov \cite{B13} in Example 1, devised a polynomial map called the \textit{additive trap} $F_{at}:\mathbb{A}_\mathbb{Z}^2\to \mathbb{A}_\mathbb{Z}^2$ by defining $F_{at}(x,y)=(x^2y,x^2y+xy^2)$. This polynomial map satisfies the following properties:
\begin{enumerate}[label=(\alph*)]
    \item It and its reductions modulo $p$ are dominant for all primes $p$.
    
    \item The only fixed point of it and any reduction of it modulo an arbitrary prime $p$ is (0,0).
    
    \item We have $F_{at}^{(p)}(x,y)\equiv (0,0)$ (mod $p$) for every $(x,y)\in \mathbb{A}^2_{\mathbb{F}_p}$ and for all primes $p$, where $F_{at}^{(p)}$ is the $p$-th iteration of $F_{at}$.
\end{enumerate}
Let $p$ be any prime. Note that all points $(x,y)\in\mathbb{A}_{\mathbb{F}_p}^2$ with either $x=0$ or $y=0$ are taken to (0,0) by $F_{at}$.  Let $x\in \mathbb{F}_p^*$. Then for any $y\in\mathbb{F}_p^*$ we get $$\frac{x^2y+xy^2}{x^2y}=\frac{y}{x}+1.$$ So after at most $p-1$ iterations the second coordinate becomes 0 and thus applying $F_{at}$ once more we reach (0,0). Since $p$ is arbitrary, we get (c). For the proofs of (a) and (b) and more details see \cite{B13}. \\
Upon further analysis of the discussion above, we notice that the $p$-th iteration of $F_{at}$ modulo $p$ is the zero map, which follows from the fact that the polynomial $u(x)=x+1$ has the following property: for every $n\in \mathbb{N}$, $u^{(n)}(x)=x+n$, so that, in particular, for every prime $p$, $u^{(p-1)}(1)=p\equiv 0$ (mod $p$). Throughout this paper by $u^{(n)}$ we will mean $u\circ u\circ \cdots\circ u$, the $n-th$ iteration of $u$. We can write our first definition now motivated from the polynomial $u(x)=x+1$ above. Suppose that $u:=u(x)\in\mathbb{Z}[x]$ is of degree $d$, $r\in\mathbb{Z}$ and $A$ is a finite subset of the set of all prime numbers. If, for every prime in $p$ not contained in $A$, there exists an $m_p\in\mathbb{N}$ such that $u^{(m_p)}(r)\equiv 0$ (mod $p$), then we will say that \textit{$u$ is weakly locally nilpotent at $r$ outside $A$}. The set of all weakly locally nilpotent polynomials at $r$ of degree $d$ will be denoted by $L_{r,A}^d$ (see other definitions in Section 2 below) and $L_{r,A}$ is the union of all such $L_{r,A}^d$, where the union is taken over $d\in\mathbb{N}$. When $A=\emptyset$, we say that $u$ is locally nilpotent at $r$. If $u$ is such that $u^{(n)}(r)=0,$ for some $n\in\mathbb{N}$, then we will say that $u$ is \textit{nilpotent at $r$} and the \textit{nilpotency index} is the least of such $n's$. We denote the set of all nilpotent polynomials at $r$ of degree $d$ and nilpotency index $i$ by $N_{r,i}^d$ and $N_r$ is the union of all such $N_{r,i}^d$, where the union is taken over $i,d\in\mathbb{N}$. Thus, for example, $u(x)=x+1\in L_{1,\emptyset}^1$ and $u(x)=x-1\in N_{1,1}^1$. Ideally, we would like to classify all weakly locally nilpotent polynomials of all possible degrees. The paper contains four main results:
\begin{enumerate}[label=(\arabic*)]
    \item Complete classification of all polynomials in $L_{r,\emptyset}$, when $r\in\{0,-1,1\}$. 
    
    \item Complete classification of all polynomials in $L_{1,A}^1$ and $L_{-1,A}^1$ for any given finite subset $A$ of the set of prime numbers. This can be found in \textbf{Theorem 3} of Section 5. To prove this we use the result of C. Corrales-Rodrig\'a\~nez and R. Schoof (see Theorem 1 \cite{CS97} and also page 6 in this paper).
    
    \item Complete classification of all polynomials in $S_r$, where $S_r=L_{r,\emptyset}\setminus N_r$. This can be found in \textbf{Corollaries 1} \& \textbf{2} in Section 4 and \textbf{Theorem 4} and \textbf{Corollary 4} of Section 5. Here also we use the result of C. Corrales-Rodrig\'a\~nez and R. Schoof and Theorem 5 \cite{TS13}.
\end{enumerate}

Theorem 5 \cite{TS13} mentioned in (3) says the following: Let $K$ be a number field, $\varphi_1,\ldots,\varphi_g:\mathbb{P}^1_K\to \mathbb{P}^1_K$ be rational maps of degrees at least 2. Let $A_1,\ldots,A_g$ be finite subsets $\mathbb{P}^1_K$ such that at most one set $A_i$ can contain a point that is not $\varphi_i-$ preperiodic, and such that there is at most one such point in that set $A_i$. Let $T_1,\ldots, T_g$ be finite subsets of $\mathbb{P}^1_K$ such that no $T_i$ contains any $\varphi_i-$preperiodic points. Then there is $M\in\mathbb{N}$ and a set of primes $\mathcal{P}'$ in $\mathcal{O}_K$ having positive density such that for any $i\in\{1,\ldots,g\}$, any $\gamma\in T_i$, any $\alpha\in A_i$, any $\mathfrak{p}\in\mathcal{P}'$ and any $m\ge M$, one has $$\varphi_i^{(m)}(\gamma)\not\equiv \alpha  (\text{mod } \mathfrak{p}).$$ Now what does this say about the locally nilpotent polynomials which are non-nilpotent at $r$? That motivates our Fact 1, which is proven in section 3.

\begin{fact}
    The $S_r$ can contain only linear polynomials.
\end{fact}

The main tools that we have used here are (see section 3 for details):
\begin{enumerate}[label=(\arabic*)]
    \item Fact 1,

    \item Fact 2,
    
    \item Lemma 1, which is a consequence of the aforementioned theorem by C. Corrales-Rodrig\'a\~nez and R. Schoof and
    
    \item The reduction of polynomials.
\end{enumerate}

This paper has 6 sections in total. Section 1 is the introduction. In sections 2 and 3 we formalize the definitions and introduce the main tools, respectively. Section 4 contains the main result listed in (1) above. Section 5 is dedicated to the classification of polynomials in $S_r$, which can only be linear polynomials by Fact 1. The last section has some open questions that arises from the study of the polynomials in this paper.

\section{Terminology and Definitions}
We will start by formally defining the polynomials mentioned in the introduction and fixing some basic terminology that we will use throughout this paper. Let $\mathcal{P}$ be the set of all primes in $\mathbb{Z}$. For a finite subset $A$ of $\mathcal{P}$ and for $a\in \mathbb{Z}$, we define
$$\mathcal{P}_A:=\mathcal{P}\setminus A\;\;\textup{and}\;\; P_A(a):=\{p\in \mathcal{P}_A~|~p \text{ divides }a\}\;\; \textup{and}\;\; P(a):=P_{\emptyset}(a).$$
So $P(a)$ is the set of all primes that divides $a$.

For $u=u(x)\in\mathbb{Z}[x]$ of degree at least 1, we define the polynomials $u^{(1)}(x):=u(x)$ and $u^{(n+1)}(x):=u(u^{(n)}(x))$, $n\in\mathbb{N}$. Having fixed $r$ in $\mathbb{Z}$ and $A$ (as above) with $d\in\mathbb{N}$ a degree and $i\in\mathbb{N}$, an index, we define the following:

\begin{enumerate}[label=(\arabic*)]

    \item We will say that $u(x)$ is a \textbf{\textit{ weakly locally nilpotent polynomial}} at $r$ outside $A$ if for each $p\in \mathcal{P}_A$, there exists $m\in\mathbb{N}$ (possibly depending on $p$) such that $u^{(m)}(r)\equiv 0$ (mod $p$). For each $p\in  \mathcal{P}_A$, we let $m_p$ be the least of all such $m'$s. We fix the following notation for weakly locally nilpotent polynomials at $r$ outside $A$:\\ 
    $L_{r,A}^d:=\{u=u(x)\in\mathbb{Z}[x]~|~u\textup{ of degree } d\textup{ is weakly locally nilpotent at }r$\\
    $\textup{ outside }A\},$\\
    $L_{r,A}:=\sqcup_{d=1}^\infty L_{r,A}^d~$.
    
    \item If $A=\emptyset$ in (1), then we will just drop the terms ``weakly" and ``outside $A$". 
    
    \item We will say that $u(x)\neq 0$ is a \textbf{\textit{nilpotent polynomial}} at $r$ if $\exists~ n\in\mathbb{N}$ such that $u^{(n)}(r)=0$. We will call the smallest of all such $n'$s as the \textbf{\textit{nilpotency index/index of nilpotency}} of $u(x)$ at $r$. If $u^{(n)}(r)\neq 0$ for all $n\in\mathbb{N}$, we will say that $u$ is \textbf{\textit{non-nilpotent}} at $r$. We fix the following notations for nilpotent polynomials at $r$:\\ $N_{r,i}^d:=\{u\in\mathbb{Z}[x]~|~u\textup{ is nilpotent at }r \textup{ of nilpotency index }i\textup{ and degree }d\},$\\
    $N_{r,i}:=\sqcup_{d=1}^\infty N_{r,i}^d~,$\\
    $N_r:=\sqcup_{i=1}^\infty N_{r,i}~.$
    
    \item The rest of the notation that we will be using are as follows:\\
    $S_r:=L_r\setminus N_r$.
    
    For integers $a,b,c\in\mathbb{Z}, \textup{ with }c\neq 0$, we will write $a\equiv_c b$ to mean $a\equiv b$ (mod $c$).
\end{enumerate}

\begin{remark}
It is clear that $N_r\subset L_{r,\emptyset}$. But it turns out that, for every given $r\in\mathbb{Z}$, $S_r$ is always non-empty (see \textbf{Corollaries 1,2,4}  and \textbf{Theorem 4} below).

\end{remark}
 
\subsection{Some examples}
\begin{enumerate}[label=(\alph*)]
    \item Let $r\in\mathbb{Z}$. For each non-zero $q(x)\in\mathbb{Z}[x]$, $(x-r)q(x)\in N_{r,1}$.
    
    \item If $u(x)=-2x-4$, then $u(-1)=-2,~u(-2)=0$. So $u(x)\in N_{-1,2}^1$. If $r\in\mathbb{Z}\setminus\{-1\}$, $u_r(x):=-(r+1)x+(r+1)^2\in N_{r,2}^1$ and if $r\in\mathbb{Z}\setminus\{0\}$, $u_r(x):=-2x+4r\in N_{r,2}^1$. Also if $r\in\mathbb{N}, u(x)=x-1\in N_{r,r}^1$.

    \item Let $u(x)=-2x^2+7x-3$. Then $u(1)=2,~u(2)=3$ and $u(3)=0$. So, $u(x)\in N_{1,3}^2$. From this and \textit{Fact 2} (stated and proved below) it follows that $v(x):=2x^2+7x+3\in N_{-1,3}^2$.

    \item The polynomial $u(x)=-x^3+9x^2-25x+25\in N_{2,4}^3$.
    
    \item \textit{This example shows the existence of non-nilpotent, locally nilpotent polynomials at 1.} Let $u(x)=x+1$. Then, by induction we see that $u^{(n)}(1)=n+1$, for every $n\in\mathbb{N}$ and hence $u\notin N_{1}$. For each $p\in \mathcal{P}$, $u^{(p-1)}(1)=p\equiv_p 0$. Thus $u(x)\in S_1$. In \textit{Corollary }1 we will see that $S_1=\{x+1\}$.
    
    \item For every $a\in\mathbb{Z}\setminus\{0\}$, let $u_a=u_a(x):=x+a$. By induction, we get $u_a^{(n)}(0)=na$. So it is clear that $u_a\notin N_0$. For each prime $p$, $u_a^{(p)}(0)=pa\equiv_p 0$. Thus $u_a\in S_0$.
    
    \item Let $u(x)=4x-2$. Then $u(1)=2$ and $u(2)=6\equiv_5 1$. This means that $u^{(n)}(1)$ is either 1 or 2 \textit{modulo} 5, for every $n\in\mathbb{N}$. This shows that $u(x)\notin L_{1,A}$, for every finite subset $A\subset \mathcal{P}_{\{5\}}$.
   
\item Let $u(x)$ be as in \textit{example }(f). Then, by induction, we have $u^{(n)}(0)=\frac{2}{3}(1-4^n)$,
which cannot be zero for any $n\in \mathbb{N}$ and so the above polynomial is not contained in $N_0$. Note that $m_2=1, ~m_3=3$. For every prime $p\in P\setminus\{2,3\}$, we have that $u^{(p-1)}(0)\equiv_p 0$ by \textit{Fermat's little theorem} and so $u(x)\in S_0$.
 \end{enumerate}

\begin{remark}
Computation of polynomial iterations is very complicated. But the linear polynomials has a nice and easy to understand iteration formula. Let $u(x)=ax+b$ be a linear polynomial, i.e., $a\in\mathbb{Z}\setminus\{0\}$. Then,
by induction, it follows that for every $n\geq 1$, $$u^{(n)}(x)=a^nx+b\left(\sum\limits_{i=0}^{n-1}a^i\right),~~n\in\mathbb{N}.$$ So, $u^{(n)}(r)=a^nr+b\left(\sum\limits_{i=0}^{n-1}a^i\right)$, for each $n\in\mathbb{N}$. Throughout this paper, we will refer to
this formula as the \textit{\textbf{linear iteration formula.}}
\end{remark} 

\section{The main tools}
\textit{In this section we will develop the necessary tools. We begin with the proof of Fact 1 which was introduced in the introduction.}

\textit{Proof of Fact 1.} Suppose that $u$ is not nilpotent at $r$ and of degree at least 2. If we take $K=\mathbb{Q}$, $g=1$, $A=A_1=\{0\}$, $T=T_1=\{r\}$ and $\varphi=\varphi_1=u$ in Theorem 5 \cite{TS13}, it follows that $u$ cannot be locally nilpotent at $r$. In other words, only the linear polynomials can be locally nilpotent without being nilpotent at $r$. So, $S_r$ can only contain linear polynomials.\endproof{}

The next fact indicates that it is enough to study the locally nilpotent polynomials at non-negative $r'$s. In particular, it shows that there is a one-to-one correspondence between $S_r$ and $S_{-r}$.

\begin{fact}
Let $u(x)$ be a polynomial of degree $d$ and let $r\in\mathbb{Z}\setminus\{0\}$. Define $v(x):=-u(-x)$. Then $u(x)\in L_{r,\emptyset}^d \iff v(x)\in L_{-r,\emptyset}^d$. Similarly $u(x)\in N_{r,n}^d \iff v(x)\in N_{-r,n}^d$ and $u(x)\in S_r\iff v(x)\in S_{-r}$.
\end{fact}

\textit{Proof.} Since $v(-x)=-u(x)$, by induction it follows that $v^{(n)}(-r)=-u^{(n)}(r)$, from which the fact follows.\endproof{}

\vspace{3mm}
Before moving on to the other tools it is imperative that we formally introduce the result by C. Corrales-Rodrig\'a\~nez and R. Schoof (\textit{Theorem }1 \cite{CS97}): Let $K$ be a number field and $x,y\in K^*$. If, for almost all prime ideals $\mathfrak{p}$ of $\mathcal{O}_K$, we have
$$\{n\in\mathbb{N}~|~y^n\equiv 1~(\textup{mod }\mathfrak{p})\}\supseteq \{n\in\mathbb{N}~|~x^n\equiv 1~(\textup{mod }\mathfrak{p})\},$$ then $\exists~ m\in\mathbb{Z}$ such that $y=x^m$. 

Now we will state and prove \textbf{Lemma 1} which was mentioned in the introduction and we will also justify its importance in understanding linear locally nilpotent polynomials.

\begin{lemma}
Let $\alpha,\beta,\gamma\in\mathbb{Z}\setminus \{0\}$ be such that neither $\frac{\beta}{\gamma}$ nor $\frac{\gamma}{\beta}$ is a non-negative power of $\alpha$. Then $\mathcal{P}\setminus \cup_{n\in\mathbb{N}}P(\gamma \alpha^n-\beta)$ is an infinite set.
\end{lemma} 

\textit{Proof.} Suppose, if possible, that $\mathcal{P}\setminus \cup_{n\in\mathbb{N}}P(\gamma \alpha^n-\beta)$ is a finite set. This means that the set $\cup_{n\in\mathbb{N}}P(\gamma \alpha^n-\beta)$ contains all but finitely many primes. Then $\cup_{n\in\mathbb{N}}P(\gamma \alpha^n-\beta)\setminus P(\gamma)$ also contains all but finitely many primes. So, for almost all $p\in \mathcal{P}_{P(\gamma)}, ~\alpha^{n_p}\equiv_p \beta\gamma^{-1}$ for some $n_p\in \mathbb{N}$ (choice of $n_p$ possibly depends on $p$). So, if $k\in\mathbb{N}$ is such that $\alpha^k\equiv_p 1$, then $(\beta \gamma^{-1})^k\equiv_p (\alpha^{n_p})^k\equiv_p 1$. Taking $\alpha=x$, $\beta\gamma^{-1}=y$ in Theorem 1 [5] applied over $\mathbb{Q}$, we arrive at a contradiction! Thus $\mathcal{P}\setminus \cup_{n\in\mathbb{N}}P(\gamma \alpha^n-\beta)$ is an infinite set.  \endproof{}

\begin{remark}[Importance of the Lemma]
Let $r\in \mathbb{Z}\setminus\{0\}$ and $u=u(x)=ax+b\in L_{r,\emptyset}^1$ with $a\neq \pm 1$. By the \textit{linear iteration formula}, we have
$$u^{(n)}(r)=\frac{a^n(r-ar-b)+b}{1-a}.$$ Since $u\in L_{r,\emptyset}^1$, we can say that $\mathcal{P}\setminus \cup_{n\in\mathbb{N}} P(\gamma \alpha^n-\beta)$ is a finite set (in fact it is an empty set), where $\alpha=a,\beta=-b$ and $\gamma=r-ar-b$. Then it follows from the above \textit{lemma} that either $\frac{\beta}{\gamma}$ or $\frac{\gamma}{\beta}$ is a power of $\alpha$, i.e., $b=-a^m(r-ar-b),$ for some $m\in\mathbb{Z}$. Moreover, if $m\in\mathbb{N}$ we can say that $u\in N_r$. So, to summarize, if $u$ is in $S_r$ (with $a\notin\{ \pm 1\}$), then $\exists~m\in\mathbb{N}\cup \{0\}$ such that $a^mb=b+ar-r$.
\end{remark}

\textbf{Reduction of polynomials.} Let $r\in \mathbb{N}$ and $u(x)\in\mathbb{Z}[x]$ such that $r|u(0)$. Define $v(x):=\frac{1}{r}u(rx)$. Note that $v$ is indeed a polynomial over $\mathbb{Z}$ of the same degree as $u(x)$ and $rv(1)=u(r)$. Using induction, one can show that $rv^{(n)}(1)=u^{(n)}(r),\forall ~n\in\mathbb{N}$. Then it follows that $u(x)$ is weakly locally nilpotent at $r$ outside $A$ iff $v(x)$ is weakly locally nilpotent at 1 outside $A\cup P(r)$ and also $u(x)$ is nilpotent at $r$ iff $v(x)$ is nilpotent at 1. Thus we can reduce any polynomial $u(x)$ in $L_{r,\emptyset}^d$ with $r|u(0)$ to the polynomial $v(x)$ in $L_{1,P(r)}^d$. We will call this the \textit{reduction of $u(x)$ to $v(x)$}.

In the following section we classify all the polynomials in $L_r$, for $r\in\{0,1,-1\}$.
\section{Arbitrary $d$ and $r\in\{0,1,-1\}$}

We will now state and prove our first main result.

\begin{theorem}
The following is the list of all polynomials in $L_{1,\emptyset}$\textup{:}
\begin{enumerate}[label=(\arabic*)]
    \item $ (x-1)p(x)$ with $p(x)\in\mathbb{Z}[x]\setminus\{0\}$ (\textit{Nilpotent of index 1}).
    
    \item $ -2x+4+p(x)(x-1)(x-2)$, with $p(x)\in\mathbb{Z}[x]$ (\textit{Nilpotent of index 2}).
    
    \item $ -2x^2+7x-3+p(x)(x-1)(x-2)(x-3)$, with $p(x)\in\mathbb{Z}[x]$ (\textit{Nilpotent of index 3}).
    
    \item $ x+1$ (\textit{Locally nilpotent but not nilpotent}).
\end{enumerate}
\end{theorem}
\textit{Proof.} Let $u=u(x)\in L_{1,\emptyset}^d$. We will consider the following three cases:

\subsubsection*{Case 1. $u(1)-1\not \in\{\pm 1\}$.}
Then $P(u(1)-1)\neq\emptyset$ and for each $p\in P(u(1)-1)$ we have $u(1)\equiv_p 1$, i.e., $m_p$ does not exist, a contradiction to that $u\in L_{1,\emptyset}^d$!

\subsubsection*{Case 2. $u(1)-1=-1$.}
These are just the polynomials listed in (1).

\subsubsection*{Case 3. $u(1)-1=1$.}
This means $u(1)=2$. Now we will explore the possibilities for $u(2)$. If $u(2)=0$, then $u(x)$ is of the form listed in (2). So suppose that $u(2)\neq 0$. Of course $u(2)\notin\{1,2\}$ as otherwise we get $u^{(n)}(1)=$ 1 or 2, for every $n\in\mathbb{N}$ and hence it cannot be in $L_{1,\emptyset}^d$. Thus $u(2)$ is either $\le -1$ or $\ge 3$, i.e., $|u(2)-1|\ge 2$. In other words, $P(u(2)-1)\neq \emptyset$. Let $p\in P(u(2)-1)$. Then $u(2)\equiv_p 1$. As $u$ is locally nilpotent at $1$, $p$ must be 2 and so $u(2)-1$ must be of the form $\pm 2^t$ for some $t\in\mathbb{N}$. To arrive at a contradiction suppose that $u(2)\neq 3$. That means $u(2)$ is either $\geq 4$ or $\leq -1$. Let's consider these two possibilities one by one.

\textbf{Possibility 1}. $u(2)\ge 4$. We know that $u(2)-1=2^t$ and so $u(2)$ is odd. So, in fact $u(2)\ge 5$. Then there exists $p\in \mathcal{P}_{\{2\}}$ such that $p\in P(u(2)-2)$. Hence  $u^{(n)}(2)\equiv_p 2$, for every $n\in\mathbb{N}$, a contradiction to the fact that $u\in L_{1,\emptyset}^d$!

\textbf{Possibility 2}. $u(2)\le -1$. We know that $u(2)-1=-2^t$ and so $u(2)$ is odd, which implies that $u(2)-2$ is odd as well and less or equal to $-3$. Using the same argument as in possibility 1 we get a contradiction!

So $u(2)$ must be 3. Next we look at $u(3)$. If $u(3)=0$, then $u(x)$ is of the form listed in (3). So suppose that $u(3)\neq 0$. For the same reason as above $u(3)\notin\{0,1,2,3\}$. Thus $u(3)$ is either $\le -1$ or $\ge 4$. To get to a contradiction suppose that $u(3)\neq 4$. Then either $u(3)-3\le -4$ or $\geq 2$. In any case, $P(u(3)-3)\neq\emptyset$. Let $p\in P(u(3)-3)$. Then $u(3)\equiv_p 3$ and so $p\in \{2,3\}$. If $p=2$ then $u(1)\equiv_p 0$. Since $3\equiv_p 1$, we must have $u(3)\equiv_p u(1)$ and so $u(3)\equiv_p 3\equiv_p 1\not\equiv_p u(1)$, which is an impossibility! So $p=3$ and $u(3)-3=\pm 3^s$, for some $s\in\mathbb{N}$. 

Again for similar reasoning as above, $P(u(3)-1)\neq \emptyset$. For each $p\in P(u(3)-1)$, we have $u(3)\equiv_p 1$ which implies $p\in\{2,3\}$. But $p|u(3)-1=2\pm 3^s$ and so $p$ cannot be 2 or 3, which is absurd! So $u(3)=4$.

Next we look at $u(4)$. We claim that no further iteration of $u$ at $1$ can be zero and we would like to prove this by showing that $u(n-1)=n$, $\forall ~n\ge 4$ and that would mean $u(x)=x+1$. We want to use \textit{mathematical induction} to prove this claim. Let $u(q-1)=q$, for every $2\le q\le n$, for some $n\ge 4$ and we want to show that $u(n)=n+1$. Since $u(1)=2, ~u(2)=3, ~u(3)=4,~\ldots, u(n-1)=n$, there is a polynomial $p(x)$ such that $u(x)=x+1+p(x)(x-1)(x-2)(x-3)\cdots (x-n+1)$. So $u(n)=n+1+p(n)\cdot (n-1)!\neq 0$, as $n\ge 4$. If $u(n)=i$ for some $i\in\{1,\ldots,n\}$, then the iterations $u^{(m)}(1)\in\{1,\ldots,n\}$, for every $m\in\mathbb{N}$ and that means for only finitely many primes $p$, $m_p$ exists. Thus $u$ cannot be locally nilpotent at $1$ and $u(n)\not \in \{0,\ldots,n\}$. This means $u(n)$ is either $\ge n+1$ or $\le -1$. For a contradiction, suppose that $u(n)\neq n+1$. Then, either $u(n)-n\ge 2$ or $u(n)-n\le -(n+1)$. In any case, we get $P(u(n)-n)\neq\emptyset$. For each $p\in P(u(n)-n)$, we have $u(n)\equiv_p n$ which is an impossibility unless $p\le n$. Suppose, if possible, $p<n$. So $n\equiv_p a$ for some $a\in \{1,\ldots,p-1\}$. Note that $a$ cannot be zero as otherwise $u(n)\equiv_p n\equiv_p 0$ and also $u(n-1)=n\equiv_p 0$. This means $p|u(0)=a_0$ and so $p|u(p)=p+1$, an impossibility! But by the induction hypothesis we have $u(a)=a+1$ and also $a\equiv_p n\equiv_p u(n)\equiv_p u(a)$, i.e., $u(a)\equiv_p a$, which is absurd as this means $m_p$ does not exist! So $p=n$, i.e., $n$ is prime and $u(n)=n\pm n^s$, for some $s\in\mathbb{N}$.

Again, for similar reasoning as above, $P(u(n)-1)\neq \emptyset$. So for every $q\in P(u(n)-1)$, $u(n)\equiv_q 1$ and it follows that $q$ is less than or equal to $n$. But if $q=n$, then $n=q|u(n)-1=(n-1)\pm n^s$ and so $n|1$, an impossibility! So, in fact, we have $q\le n-1$. We can choose $b\in\{0,\ldots,q-1\}$ such that $n\equiv_q b+1$. By the induction hypothesis $u(b+1)=b+2$ and also $u(b+1)\equiv_q u(n)\equiv_q 1$. These two relations together imply $b+1\equiv_q 0$, i.e., $q|n$. But, since $n$ is a prime, $n=q$, which is a absurd as $q\le n-1$. Thus $u(n)=n+1$.\endproof{}

\begin{remark}
It follows from \textbf{Fact 2} and \textbf{Theorem 1} that the following polynomials are in $L_{-1,\emptyset}$:
\begin{enumerate}[label=(\arabic*)]
    \item $ (x+1)p(x)$, with $p(x)\in\mathbb{Z}[x]\setminus\{0\}$.
    
    \item $ -2x-4+p(x)(x+1)(x+2)$, with $p(x)\in\mathbb{Z}[x]$.
    
    \item $ 2x^2+7x+3+p(x)(x+1)(x+2)(x+3)$, with $p(x)\in\mathbb{Z}[x]$.
    
    \item $ x-1$.
\end{enumerate}
\end{remark}

\begin{corollary}
The sets $S_1$ and $S_{-1}$ are singleton sets.
\end{corollary} 

\textit{Proof.} Let $u(x)\in S_1$. Then by \textbf{Theorem 1}, $u(x)$ must be $x+1$ as all the other polynomials in the list (1)-(4) of \textbf{Theorem 1} are in $N_1$. Now by \textit{Fact 2}, it follows that $S_{-1}=\{x-1\}$. \endproof{}

\begin{theorem}
The following is the list of all polynomials in $L_{0,\emptyset}$\textup{:}
\begin{enumerate}[label=(\arabic*)]
    \item $ ax$, with $a\in\mathbb{Z}\setminus\{0\}$ (\textit{Locally nilpotent but not nilpotent}).
    
    \item $ \pm x+b$, with $b\in\mathbb{Z}\setminus \{0\}$ (\textit{Locally nilpotent but not nilpotent}).
    
    \item $ ax+b$, with $\mathcal{P}\supsetneq P(b)\supseteq P(a)\neq \emptyset$ (\textit{Locally nilpotent but not nilpotent}). 

    \item $xp(x)$, with $p(x)\in\mathbb{Z}[x]\setminus\{0\}$ (\textit{Nilpotent of index 1}).

    \item $(x-a)p(x)$, with $a\in\mathbb{Z}\setminus\{0\}$ and $p(x)\in \mathbb{Z}[x]$ s.t $p(0)=-1$ (\textit{Nilpotent of index 2}).
    
\end{enumerate}
\end{theorem}

\textit{Proof.} First suppose that $u$ is nilpotent of nilpotency index $m$, for some $m\in\mathbb{N}$. If $u(0)=0$, then $m=1$ and $u(x)=xp(x),$ for some non-zero $p(x)\in\mathbb{Z}[x]$, which is (4) in the list. So, suppose that $u(0)\neq 0$. Define $$u_0:=u(0),~u_n:=u^{(n+1)}(0)-u^{(n)}(0),~n\in\mathbb{N}.$$ Then $u_{n+1}=u^{(n+2)}(0)-u^{(n+1)}(0)= u(u^{(n+1)}(0))-u(u^{(n)}(0)).$ That means $u_n$ divides $u_{n+1},\forall~n\in\mathbb{Z}_{\ge 0}$. We also have $u^{(m)}(0)=0$ and so $u_m=u^{(m+1)}(0)-u^{(m)}(0)=u^{(m+1)}(0)=u_0$. 
As $u_0|u_1|\ldots|u_m=u_0,$ it follows that $u_n=\pm u_0,$ for all $n.$ Note that $u_0+\cdots+u_{m-1}=u^{(m)}(0)=0$. This means $m$ must be even and half these integers are positive and the other half are negative (since $|u_n|=|u_0|,\forall~n\in\mathbb{N}$). So there exists $k\in\{1,\ldots,m-1\}$ such that $u_{k-1}=-u_k,$ i.e., $u^{(k)}(0)-u^{(k-1)}(0)=u^{(k)}(0)-u^{(k+1)}(0),$ i.e., $u^{(k+1)}(0)=u^{(k-1)}(0)$. Thus $u^{(n+2)}(0)=u^{(n)}(0),\forall ~n\ge k-1$ and so in particular, we have $0=u^{(m)}(0)=u^{(m+2)}(0)=u^{(2)}(0)$. Hence, $m=2$ and $u(x)=(x-\alpha)p(x)$, with $\alpha\in\mathbb{Z}\setminus\{0\}$ and $p(x)\in\mathbb{Z}[x]$ with $p(0)=-1$; here $\alpha=u(1),$ which is (5) in the list.

Now suppose that $u\in S_0$. Then by Fact 1, it must be linear. Let $u(x):=ax+b, a\neq 0$. Note that if $b=0$ we get (1) listed above. Now suppose $b\neq 0$, i.e., $P(b)\neq \mathcal{P}$. When $a=1$, every $u(x)\in S_0$: in fact if $u(x)=x+b$, then by the \textit{linear iteration formula}, $u^{(n)}(0)=b(1+\cdots+1)=bn$, which is always non-zero, $\forall~n\in\mathbb{N}$ and for each prime $p$, $u^{(p)}(0)=bp\equiv_p 0$. When $a=-1$, $u^{(2)}(0)=0$. So $-x+b\in N_{0,2}^1$.\\
So we can assume that $|a|\geq 2$, i.e., $P(a)$ is a non-empty, finite set. Again, using the \textit{linear iteration formula}, we get $u^{(n)}(0)=b(1+\cdots+a^{n-1}),~n\in\mathbb{N}$. Suppose, if possible, $p\in P_{P(b)}(a)$. Then $u^{(n)}(0)\equiv_p b$, for every $n\in\mathbb{N}$ and so $u(x)$ cannot be locally nilpotent. This means that $P(b)\supseteq P(a)\neq\emptyset$. If $p\in P(b)$, it can be checked that $m_p=1$.\\
If $p\notin P(b)\cup P(a-1)$, 
then $u^{(p-1)}(0)=\frac{b}{a-1}(a^{p-1}-1)\equiv_p 0$.\\
Finally if $p\in P(a-1)$, then $u^{(p)}(0)=b(1+\cdots+a^{p-1})\equiv_p b(1+\cdots+1)\equiv_p 0$. Thus $m_p$ exists for every $p\in\mathcal{P}$.
\endproof{}

\begin{corollary}
The following is the list of all polynomials in $S_{0,\emptyset}$\textup{:}
\begin{enumerate}[label=(\arabic*)]
    \item $ ax$, with $a\in\mathbb{Z}\setminus\{0\}$ (\textit{Locally nilpotent but not nilpotent}).
    
    \item $ \pm x+b$, with $b\in\mathbb{Z}\setminus \{0\}$ (\textit{Locally nilpotent but not nilpotent}).
    
    \item $ ax+b$, with $\mathcal{P}\supsetneq P(b)\supseteq P(a)\neq \emptyset$ (\textit{Locally nilpotent but not nilpotent}).
\end{enumerate}
\end{corollary}

\section{Linear case $d=1$}

We now state and prove our next main result:

\begin{theorem}
Let $q_1,\ldots, q_k$ be $k$ distinct primes and $A=\{q_1,\ldots, q_k\}$. Then the following is the list of all the polynomials in $L_{1,A}^1$\textup{:}
\begin{enumerate}[label=(\arabic*)]
    \item $ x\pm q_1^{s_1}\cdots q_k^{s_k}$, where $s_i\in\mathbb{N}\cup\{0\}$.
    
    \item $ \alpha(x-1)$, $\alpha\in\mathbb{Z}\setminus\{0\}$.
    
    \item $ \pm q_1^{s_1}\cdots q_k^{s_k} x+1$, , where $s_i\in\mathbb{N}\cup\{0\}$ such that $\sum s_i\ge 1$.
    
    \item $ -2x-1$ (only when $2\in A$).
    
    \item $ -2x+4$.
\end{enumerate}
\end{theorem}

\textit{Proof.} Let $u=u(x):=ax+b\in L_{1,A}^1$. It is clear that we can assume $b\neq 0$, otherwise $u^{(n)}(1)=a^n$ cannot be divisible by any prime $p\in \mathcal{P}\setminus P(a)$. By the \textit{linear iteration formula}, $u^{(n)}(1)=a^n+b(1+\cdots+a^{n-1})$, for every $n\in\mathbb{N}$. 
\begin{align*}
    \text{Note that if } a=1, ~u^{(m_p)}(1)=& 1+bm_p\equiv_p 0, \text{ for every prime }p\notin A,\\
    \implies & bm_p\equiv_p -1, \text{ for every prime }p\notin A,\\
    \implies & b\text{ is invertible in }\mathbb{F}_p, \text{ for every prime }p\notin A,\\
    \implies & b=\pm q_1^{s_1}\cdots q_k^{s_k}, \text{ for some }s_i's \textup{ in }\mathbb{N}\cup \{0\}.
\end{align*}

If $a=-1$, $u(x)=-x+b$ and $u^{(2)}(x)=x$. So $u$ cannot be in $L_{1,A}^1$  unless $b=1$ and in that case it is in fact in $N_{1,1}^1$. Thus we can assume that $|a|\ge 2$. Similar to \textbf{Theorem 1} we can break down these polynomials into the following three cases:

\subsubsection*{Case 1. $u(1)-1\notin \{\pm 1\}$.}
This means that $P(u(1)-1)\neq \emptyset$. So $a+b=1 \pm q_1^{s_1}\cdots q_k^{s_k}$, i.e., $b=1-a \pm q_1^{s_1}\cdots q_k^{s_k}$, for some $s_i\in\mathbb{N}\cup\{0\}$ with $\sum s_i\ge 1$. Then by the \textit{linear iteration formula}, we have $$u^{(n)}(1)=\frac{b\pm a^n(1-a-b)}{1-a}$$ and it follows from \textit{Remark} 3 that $\exists ~m\in\mathbb{Z}$ such that $b=\pm a^m(1-a-b)$. If $m=0$, then $b=\pm(1-a-b)$, i.e., $a+2b=1$ or $a=1$. Since $|a|\ge 2$, we deduce that $a\neq 1$ and so $a+2b=1$. Also we have $a+b=1\pm q_1^{s_1}\cdots q_k^{s_k}$. Solving $a$ and $b$ from these two equations we get $a=1\pm 2q_1^{s_1}\cdots q_k^{s_k},~b=\mp q_1^{s_1}\cdots q_k^{s_k}$. So, $u(x)=(1\pm 2q_1^{s_1}\cdots q_k^{s_k})x\mp q_1^{s_1}\cdots q_k^{s_k}$ and so $u^{(n)}(1)=\frac{1+(1\pm 2q_1^{s_1}\cdots q_k^{s_k})^n}{2},~n\in\mathbb{N}$. Letting $\alpha=1\pm 2q_1^{s_1}\cdots q_k^{s_k},\beta=-1$ and $\gamma =1$, it is clear that neither $\frac{\beta}{\gamma}$ nor $\frac{\gamma}{\beta}$ is a power of $\alpha$. Thus, it follows from \textit{Lemma }1 that $(1\pm 2q_1^{s_1}\cdots q_k^{s_k})x\mp q_1^{s_1}\cdots q_k^{s_k}\notin L_{1,A}^1$. So $|m|\ge 1$.

\underline{If $m\in\mathbb{N}$ and $b=a^m(1-a-b)$}, we get $b(1+a^m)=a^m(1-a)$. Since $\gcd(a^m,a^m+1)=1$, we must have $a^m+1|1-a$ which is only possible if $m=1$.

\underline{If $m\in\mathbb{N}$ and $b=-a^m(1-a-b)$}, we get $b(1-a^m)=-a^m(1-a)$, i.e., $b(1+\cdots+a^{m-1})=-a^m$. Since $\gcd(1+\cdots+a^{m-1},a^m)=1$, we must have $1+\cdots+a^{m-1}=\pm 1$ which is only possible if $m\in\{1,2\}$.

\underline{If $m=-n,\text{ with }n\in\mathbb{N}$ and $b=a^m(1-a-b)$}, we get $ba^n=1-a-b$, i.e., $b(a^n+1)=1-a$. It follows from above that this is only possible if $n=1$.

\underline{If $m=-n,\text{ with }n\in\mathbb{N}$ and $b=-a^m(1-a-b)$}, we get $ba^n=-(1-a-b)$, i.e., $b(a^n-1)=a-1$. Again using the same logic as above, we conclude that $n\in\{1,2\}$.

Thus we only need to look at the following four subcases:

\textit{Subcase 1.} $m=-1$.\\
Here we have $ba=\pm(1-a-b)$. First suppose that $ba=1-a-b$, i.e., $b(a+1)=1-a$. This means that $a+1|a-1$ and this is only possible if $a=-2$ and $a=-3$. These values generate the polynomials $u(x)=-2x-3$ and $u(x)=-3x-2$, respectively. When $u(x)=-2x-3$, the \textit{linear iteration formula} gives $$u^{(n)}(1)=2(-2)^n-1.$$ Letting $\alpha=-2,\beta=1$ and $\gamma=2$, it is clear that neither $\frac{\beta}{\gamma}$ nor $\frac{\gamma}{\beta}$ is a power of $\alpha$. So, by \textit{Lemma }1, $-2x-3\notin L_{1,A}^1$. Similarly we can show that $-3x-2\notin L_{1,A}^1.$

Now suppose $ba=-(1-a-b)$. This gives $b=1$ and hence $a=\pm q_1^{s_1}\cdots q_k^{s_k}$. Thus $u(x)=\pm q_1^{s_1}\cdots q_k^{s_k}x+1$ and it follows from the \textit{linear iteration formula} that $$u^{(n)}(1)=(\pm q_1^{s_1}\cdots q_k^{s_k})^n+[1+\cdots+(\pm q_1^{s_1}\cdots q_k^{s_k})^{n-1}]=\frac{1-(\pm q_1^{s_1}\cdots q_k^{s_k})^{n+1}}{1-(\pm q_1^{s_1}\cdots q_k^{s_k})}, ~n\in\mathbb{N}.$$ If $p\in P_A(1-(\pm q_1^{s_1}\cdots q_k^{s_k}))$, then $u^{(p)}(1)\equiv_p p\equiv_p 0$. So let $p\notin P_A(1-(\pm q_1^{s_1}\cdots q_k^{s_k}))$. Now, if $2\in A$, then the existence of $m_2$ is not a concern and if $2\notin A$, then $2\in P_A(1-(\pm q_1^{s_1}\cdots q_k^{s_k}))$ which was covered above.\\
Finally, if $p\notin \mathcal{P}_{A\cup \{2\}}$, then $u^{(p-2)}(1)\equiv_p 0$, by \textit{Fermat's Little Theorem}.

\textit{Subcase 2.} $m=1$.\\
Here we have $b=\pm a(1-a-b)$. First suppose that $b=a(1-a-b)$, i.e., $b(a+1)=a(1-a)$. Same reasoning as above implies $a+1|a-1$ so that only possibilities we get are $a=-2,~b=6$ or $a=-3,~b=6$. These values produces the polynomials $u(x)=-2x+6$ and $u(x)=-3x+6$, respectively. When $u(x)=-2x+6$, the \textit{linear iteration formula} gives $$u^{(n)}(1)=-(-2)^n+2.$$ Letting $\alpha=-2,\beta=-2$ and $\gamma=-1$, it is clear that neither $\frac{\beta}{\gamma}$ nor $\frac{\gamma}{\beta}$ is a power of $\alpha$. So, by \textit{Lemma }1, $-2x+6\notin L_{1,A}^1$. Similarly, we can show that $-3x+6\notin L_{1,A}^1$.

\textit{Subcase 3.} $m=-2$.\\
Here we have $ba^2=\pm (1-a-b)$. First suppose $ba^2=1-a-b$, i.e., $b(a^2+1)=1-a$. This means that $a^2+1|a-1$ which is not possible as $|1-a|\le 1+|a|<1+a^2$. Thus $ba^2=-(1-a-b)$, i.e., $b(a+1)=1$, i.e., $b=a+1=\pm 1$. So $u(x)=-2x-1$. It follows from the \textit{linear iteration formula} that $$u^{(n)}(1)=(-2)^n-[1+\cdots+(-2)^{n-1}]=\frac{(-2)^{n+2}-1}{3},~n\in\mathbb{N}.$$ It is easy to see that $m_2$ does not exist, $m_3=1$ and for all $p\in \mathcal{P}_{\{2,3\}}$, $u^{(p-3)}(1)\equiv_p 0$. So $-2x-1$ is in $L_{1,A}^1$ iff $2\in A$.

\textit{Subcase 4.} $m=2$.\\
Here we have $b=\pm a^2(1-a-b)$. First suppose that $b=a^2(1-a-b)$, i.e., $b(1+a^2)=a^2(1-a)$. Since $\gcd(1+a^2,a^2)=1$, $1+a^2|1-a$ which is impossible (see the above subcase). Now suppose that $b=-a^2(1-a-b)$, i.e., $b(a+1)=-a^2$ which means $a+1=\pm 1$ and $b=\pm a^2$. Since $|a|\ge 2$, this means $a=-2$ and $b=4$. But then $u(1)-1=1\in \{\pm 1\},$ an impossibility in this \textit{case}!

\subsubsection*{Case 2. $u(1)-1=-1$, i.e., $u(1)=0$.}
These are the polynomials in $N_{1,1}^1$.

\subsubsection*{Case 3. $u(1)-1=1$, i.e., $u(1)=2$.}
If $u(2)=0$, then $u(x)=-2x+4$. So, we can suppose that $u(2)\notin\{0,1,2\}$, i.e., $u(2)$ is either $\le -1$ or $\ge 3$, i.e., $|u(2)-1|\ge 2,$ i.e., $P(u(2)-1)\neq \emptyset$. If $u(2)=3$, then $u(x)=x+1\in S_1$. So we can further assume that $u(2)\neq 3$. Since $u(1)=2,~b=2-a$ and so $u(x)=ax+(2-a)$. Then by the \textit{linear iteration formula}, we get
$$u^{(n)}(1)=\frac{2-a-a^n}{1-a},~~n\in\mathbb{N}.$$
Since $u\in L_{1,A}^1$ it follows from \textit{Lemma }1 that $2-a=a^m$, for some $m\in\mathbb{Z}$.
If $m=0$, then $2-a=1$, i.e., $a=1$, which is not possible. Also note that if $m=-n \textup{ for some}~n\in\mathbb{N}, a^n(2-a)=1$. But this is an impossibility as $|a|\ge 2$. Thus $m\in\mathbb{N}$ and $2=a(1+ a^{m-1})$. Therefore $a=\pm 2$ and $1+a^{m-1}=\pm 1,$ i.e., $a^{m-1}=-2$, i.e., $a=-2$. But $a=-2$ implies that $b=4$ and hence $u(2)=2a+b=0$, which cannot happen as we have already said that $u(2)\not\in\{0,1,2,3\}$. Thus $ax+(2-a)\notin L_{1,A}^1$. \endproof{}

The next corollary follows directly from the computations in the proof of \textbf{Theorem 3}:

\begin{corollary}
Let $q_1,\ldots,q_k$ be $k$ distinct primes and $A=\{q_1,\ldots,q_k\}$. Then the following is the list of all polynomials in $L_{1,A}^1\setminus N_1$\textup{:}
\begin{enumerate}[label=(\arabic*)]
    \item $ x+ q_1^{s_1}\cdots q_k^{s_k}$, where $s_i\in\mathbb{N}\cup\{0\}$.
    
    \item $ x- q_1^{s_1}\cdots q_k^{s_k}$, where $s_i\in\mathbb{N}\cup\{0\}$ such that $\sum s_i\ge 1$.
    
    \item $ \pm q_1^{s_1}\cdots q_k^{s_k} x+1$, where $s_i\in\mathbb{N}\cup\{0\}$ such that $\sum s_i\ge 1$.
    
    \item $ -2x-1$ (only when $2\in A$).
\end{enumerate}
\end{corollary}

Finally we state and prove the last (main) result of this paper.

\begin{theorem}
Let $r=q_1^{a_1}\cdots q_k^{a_k}$ be the prime decomposition of $r$. Then the following is the list of all polynomials in $S_r$\textup{:}
\begin{enumerate}[label=(\arabic*)]
    \item $x+ q_1^{s_1}\cdots q_k^{s_k}$, where $s_i\in\mathbb{N}\cup\{0\}$.
    
    \item $x-q_1^{s_1}\cdots q_k^{s_k}$, where $s_i\in\mathbb{N}\cup\{0\}$ with at least one $j\in\{1,\ldots, k\}$ s.t. $s_j>a_j$.
    
    \item $\pm q_1^{s_1}\cdots q_k^{s_k}x+r$, where $s_i\in\mathbb{N}\cup\{0\}$ with $\sum\limits_{i} s_i\ge 1$.
    
    \item $-2x-r$ (only when $r$ is even).
\end{enumerate}
\end{theorem}

\textit{Proof.} Let $u=u(x):=ax+b\in S_r$ and $A:=P(r)$. First we will look at the instances when $a=\pm 1$.
\begin{align*}
    \text{Note that if } a=1, ~u^{(m_p)}(1)=&1+bm_p\equiv_p 0, \text{ for every prime }p\notin \{q_1,\ldots, q_k\},\\
    \implies & bm_p\equiv_p -1, \text{ for every prime }p\notin \{q_1,\ldots, q_k\},\\
    \implies & b\text{ is invertible in }\mathbb{F}_p, \text{ for every prime }p\notin \{q_1,\ldots, q_k\},\\
    \implies & b=\pm q_1^{s_1}\cdots q_k^{s_k}, \text{ for some }s_i's \textup{ in }\mathbb{N}\cup \{0\}.
\end{align*}
So, $u(x)$ is the form $x \pm q_1^{s_1}\cdots q_k^{s_k}$. First suppose that $u(x)=x+q_1^{s_1}\cdots q_k^{s_k}$. Then by the \textit{linear iteration formula}, $u^{(n)}(r)=q_1^{a_1}\cdots q_k^{a_k}+n\cdot q_1^{s_1}\cdots q_k^{s_k}$, which is always non-zero for every $n\in\mathbb{N}$, which are the polynomials in (1) in the list above. Now suppose that $u(x)=x-q_1^{s_1}\cdots q_k^{s_k}$. If, for all $i\in\{1,\ldots,k\}$, $s_i\le a_i$, then $u^{(q_1^{a_1-s_1}\cdots q_k^{a_k-s_k})}(r)=0$, which is a contradiction as $u$ is non-nilpotent! That means we must have at least one $j\in\{1,\ldots,k\}$ such that $a_j<s_j$. Then it can be easily checked that $u^{(n)}(r)$ can never be zero for any $n\in\mathbb{N}$ and so we get the polynomials in (2) in the list above.

If $a=-1$, $u(x)=-x+b$ and $u^{(2)}(x)=x$. So, $u$ cannot be in $L_{r,\emptyset}^1$  unless $b=r$ and in that case it is in fact in $N_{r,1}^1$, a contradiction! So $a\neq -1$. Thus we can $|a|\ge 2$. It follows from \textit{Remark }3 that $\exists~m\in\mathbb{N}\cup \{0\}$ such that $a^mb=b+ar-r$. If $m=0$ then $r(1-a)=0$ which is impossibility as $r\neq 0$ and $|a|\ge 2$. That means that $m\in\mathbb{N}$. Thus $b(a^m-1)=r(a-1),$ i.e., $b(1+\cdots+a^{m-1})=r$. This means $b|r$. 

We now want to show that $u(r)-r\notin \{\pm 1\}$. Suppose not. Then $u(r)=r\pm 1$. This means that $b=r-ar\pm 1$, i.e., $u(x)=ax+(r-ar\pm 1)$. We will only consider the possibility that $b=r-ar-1$ as the other possibility can be rejected using the same argument. Applying the \textit{linear iteration formula}, we get $$u^{(n)}(r)=\frac{a^n+r-ar-1}{1-a}=\frac{a^n+b}{1-a},~\text{ for all }n\in \mathbb{N}.$$ From \textit{Remark} 3 it follows that $r-ar-1=-a^m$, for some $m\in\mathbb{Z}$. It is clear that $m\neq 0$, as otherwise $r-ar=0,$ i.e., $r(1-a)=0,$ i.e., either $r=0$ or $a=1$, which is not true!
If $m=-n$ for some $n\in\mathbb{N}$, then $a^n(r-ar-1)=-1,$ again an impossibility as $|a|\ge 2$! Thus $m\in\mathbb{N}$ and $r(1-a)=1-a^m$, i.e., $r=1+\cdots+a^{m-1}$. So $m\ge 2$, $a|r-1$ and $u^{(m)}(r)=0$, i.e., $u$ is nilpotent at $r$, a contradiction! This proves that $u(r)-r$ cannot be a unit. That means $u(r)=r\pm q_1^{s_1}\cdots q_k^{s_k}$, for a suitable collection of $s_i's$ in $\mathbb{N}\cup\{0\}$ with $\sum\limits_i s_i\ge 1$. So $b=r-ar\pm q_1^{s_1}\cdots q_k^{s_k}$. But then $b|r$ implies that $b|q_1^{s_1}\cdots q_k^{s_k}$, i.e., $\exists~t_i\in\mathbb{N}\cup \{0\}$ with $s_i\ge t_i$ for every $i$ such that $b=\pm q_1^{t_1}\cdots q_k^{t_k}$. From $b=r-ar\pm q_1^{s_1}\cdots q_k^{s_k}$ we get $ra=r-b\pm q_1^{s_1}\cdots q_k^{s_k}=r-b(\pm 1\pm q_1^{s_1-t_1}\cdots q_k^{s_k-t_k})$, i.e., $r|b(\pm 1\pm q_1^{s_1-t_1}\cdots q_k^{s_k-t_k})$.\\
Suppose, if possible, all the $t_i's$ are zero. Then $b=\pm 1$ and so $r$ must divide $\pm 1\pm q_1^{s_1}\cdots q_k^{s_k}$ but this is clearly absurd as $\gcd(r,\pm 1\pm q_1^{s_1}\cdots q_k^{s_k})=\gcd(q_1^{a_1}\cdots q_k^{a_k},\pm 1\pm q_1^{s_1}\cdots q_k^{s_k})=1$. So $\sum\limits_i t_i\ge 1$. All this now boils down to the following two cases:

\subsubsection*{Case 1. $\exists ~j\in\{1,\ldots,k\}$ s.t. $s_j>t_j$.}
Since $\gcd(r,\pm 1\pm q_1^{s_1-t_1}\cdots q_k^{s_k-t_k})=1$, $r|b$, i.e., $r=\pm b$ (since we already had $b|r$). So $a_i=t_i\le s_i,~\forall~i\in\{1,\ldots,k\}$. So we can use the reduction of polynomials. Define $$v=v(x):=\frac{1}{r}u(rx)=ax\pm 1.$$ Then $v(1)=1\pm q_1^{s_1-a_1}\cdots q_k^{s_k-a_k}$ and $v\in L_{1,A}^1\setminus N_1$. It follows from the list in \textit{Corollary }3 that we have $2$ possibilities for $v$:
\begin{enumerate}[label=(\roman*)]
    
    \item $v(x)=\pm q_1^{s_1-a_1}\cdots q_k^{s_k-a_k}x+1$. Then $u(x)=\pm q_1^{s_1-a_1}\cdots q_k^{s_k-a_k}x+r$.
    
    \item $v(x)=-2x-1$ (only when $2\in A$). Then $u(x)=-2x-r$.
\end{enumerate}
One can check easily that both (i) and (ii) are indeed in $S_r$.

\subsubsection*{Case 2. $s_i=t_i,~\forall~i\in\{1,\ldots, k\}$.}
Then $\pm q_1^{s_1}\cdots q_k^{s_k}=b|r$, i.e., $a_i\ge s_i$ for each $i$. 
 From $b=r-ar\pm q_1^{s_1}\cdots q_k^{s_k}$ we get $r(1-a)=\pm 2q_1^{s_1}\cdots q_k^{s_k}=\pm 2b$. Thus either $r=\pm b$ or $r=\pm 2b$. The first possibility has been taken care of in \textit{Case }2. So we can assume that $r=2q_1^{s_1}\cdots q_k^{s_k}=\pm 2b$. But that means $2\in A$. Without loss of generality, we can assume that $q_1=2$ and so $r=2^{s_1+1}\cdots q_k^{s_k}$. Rewriting $b=r-ar\pm q_1^{s_1}\cdots q_k^{s_k}$ gives us $ra=r-2b.$ Since $ra\neq 0,$ $r=-2b$ and so $ra=-4b$, i.e., $r=-2b$ and $a=2$. This means that $u(x)=2x-\frac{r}{2}$. It follows from the \textit{linear iteration formula} that $$u^{(n)}(r)=r\cdot \frac{2^n+1}{2},~n\in\mathbb{N}.$$ Letting $\alpha=2,\beta =-1, \gamma=1$, we can see that neither $\frac{\beta}{\gamma}$ nor $\frac{\gamma}{\beta}$ is a power $\alpha$. Thus from \textit{Lemma }1 it follows that $2x-\frac{r}{2}\not\in S_r$. This completes the proof.\endproof{}

\vspace{5mm}

It follows directly from \textit{Fact 2} and \textbf{Theorem 4} that:
\begin{corollary}
If $r=-q_1^{a_1}\cdots q_k^{a_k}$ is the prime decomposition of an integer $r\le -2$, then the following is the list of all polynomials in $S_r$\textup{:}
\begin{enumerate}[label=(\arabic*)]
    \item $x- q_1^{s_1}\cdots q_k^{s_k}$, where $s_i\in\mathbb{N}\cup\{0\}$.
    
    \item $x+q_1^{s_1}\cdots q_k^{s_k}$, where $s_i\in\mathbb{N}\cup\{0\}$ with at least one $j\in\{1,\ldots, k\}$ s.t. $s_j>a_j$.
    
    \item $\pm q_1^{s_1}\cdots q_k^{s_k}x-r$, where $s_i\in\mathbb{N}\cup\{0\}$ with $\sum\limits_{i} s_i\ge 1$.
    
    \item $-2x+r$ (only when $r$ is even).
\end{enumerate}
\end{corollary}

\section{Some open problems}
For $u(x)\in\mathbb{Z}[x]\setminus\{0\}$, let
$$N(u):=\{r\in\mathbb{Z}~|~u\in N_r \} ,~~LN(u):=\{r\in\mathbb{Z}~|~u\in L_{r,\emptyset} \}.$$

\begin{enumerate}[label=Q\arabic*.]
    \item Describe all $u's$ such that $N(u)$ is finite.
    
    \item Describe all $u's$ such that $LN(u)$ is finite.
    
    \item Given $r\in\mathbb{Z}$, describe all $u's$ such that $r\in LN(u)$.
\end{enumerate}

\section{Acknowledgements} The author gratefully acknowledges his advisors Prof. Alexander Borisov and Prof. Adrian Vasiu for their constant support, their encouragements and very helpful suggestions. The author would also like to thank Prof. Jeremy Rouse for suggesting \cite{CS97} which was one of the main tools used to prove Theorems 3 and 4 and Prof. Kiran Kedlaya for maintaining a wonderful archive of William Lowell Putnam Mathematics Competition questions and answers here \url{https://kskedlaya.org/putnam-archive/}, which was very helpful in the proof of Theorem 2. Finally, the author would like to thank Prof. Thomas Tucker for his suggestion to use Theorem 5 of the paper \cite{TS13}, which led to Fact 1 and that was very useful to minimize a lot of computations.

\end{document}